\documentclass[a4paper,notitlepage,10pt]{article}%
\usepackage{amssymb}
\usepackage{graphicx}
\usepackage{amsmath}
\usepackage{amsthm}
\usepackage{amsfonts}%
\providecommand{\U}[1]{\protect\rule{.1in}{.1in}}
\theoremstyle{plain}
\newtheorem{theorem}{Theorem}[section]

\newtheorem{corollary}[theorem]{Corollary}

\newtheorem{lemma}[theorem]{Lemma}

\newtheorem{proposition}[theorem]{Proposition}

\theoremstyle{definition}

\newtheorem{remark}[theorem]{Remark}
\numberwithin{equation}{section}
\numberwithin{theorem}{section}
\ifx\pdfoutput\relax\let\pdfoutput=\undefined\fi
\newcount\msipdfoutput
\ifx\pdfoutput\undefined\else
\ifcase\pdfoutput\else
\ifx\paperwidth\undefined\else
\ifdim\paperheight=0pt\relax\else\pdfpageheight\paperheight\fi
\ifdim\paperwidth=0pt\relax\else\pdfpagewidth\paperwidth\fi
\fi\fi\fi
\begin{document}

\title{Positivity results for indefinite sublinear elliptic problems via a continuity
argument 
 }
\author{U. Kaufmann\thanks{FaMAF, Universidad Nacional de C\'{o}rdoba, (5000)
C\'{o}rdoba, Argentina. \textit{E-mail address: }kaufmann@mate.uncor.edu} , H.
Ramos Quoirin\thanks{Universidad de Santiago de Chile, Casilla 307, Correo 2,
Santiago, Chile. \textit{E-mail address: }humberto.ramos@usach.cl} , K. Umezu \thanks{Department of Mathematics, Faculty
of Education, Ibaraki University, Mito 310-8512, Japan. \textit{E-mail
address: }kenichiro.umezu.math@vc.ibaraki.ac.jp}
\and \noindent}
\maketitle

\begin{abstract}
We establish a positivity property for a class of semilinear elliptic problems
involving indefinite sublinear nonlinearities. Namely, we show that any
nontrivial nonnegative solution is positive for a class of problems the strong
maximum principle does not apply to. Our approach is based on a continuity
argument combined with variational techniques, the sub and supersolutions method and some
\textit{a priori} bounds. Both Dirichlet and Neumann homogeneous boundary
conditions are considered. As a byproduct, we deduce some existence and
uniqueness results. Finally, as an application, we derive some positivity
results for indefinite concave-convex type problems.
\end{abstract}

%\maketitle

\section{Introduction}

Let $\Omega$ be a bounded and smooth domain of $\mathbb{R}^{N}$ with $N\geq1$.
The purpose of this article is to discuss the existence of positive solutions
for the problems
\[
\left\{
\begin{array}
[c]{lll}%
-\Delta u=a(x)f(u) & \mathrm{in} & \Omega,\\
u=0 & \mathrm{on} & \partial\Omega,
\end{array}
\right.  \leqno{(P)}
\]
and
\[
\left\{
\begin{array}
[c]{lll}%
-\Delta u=a(x)f(u) & \mathrm{in} & \Omega,\\
\frac{\partial u}{\partial\nu}=0 & \mathrm{on} & \partial\Omega,
\end{array}
\right.  \leqno{(P')}
\]
where $\nu$ is the outward unit normal to $\partial\Omega$.

Here $a\in L^{r}\left(  \Omega\right)  $, $r>N$, is a function that
\textbf{changes sign}, and $f:[0,\infty)\rightarrow\lbrack0,\infty)$ is
continuous and sublinear in the following sense:
\[
\lim_{s\rightarrow0^{+}}\frac{f(s)}{s}=\infty\quad\text{and}\quad
\lim_{s\rightarrow\infty}\frac{f(s)}{s}=0.\leqno{(H_1)}
\]
In addition, we assume that $f(s)>0$ for $s>0$. The model for such $f$ is
$f(s)=s^{q}$ with $0\leq q<1$.

By a \textit{nonnegative} \textit{solution}  of
$(P)$ we mean a function $u\in W^{2,r}\left(  \Omega\right)  \cap W_{0}%
^{1,r}\left(  \Omega\right)  $ (and thus $u\in\mathcal{C}^{1}(\overline
{\Omega})$) that satisfies the equation for the weak derivatives and $u\geq0$
in $\Omega$. If, in addition, $u>0$ in $\Omega$, then we call it a
\textit{positive solution} of $(P)$. Similarly, by a \textit{nonnegative
solution} of $(P^{\prime})$ we mean a function $u\in W^{2,r}\left(
\Omega\right)  $ that satisfies the equation for the weak derivatives and the
boundary condition in the usual sense, and such that $u\geq0$ in $\Omega$. If,
in addition, $u>0$ in $\Omega$, then we call it a \textit{positive solution}
of $(P^{\prime})$.

%Let us point out that a basic difference between $(P)$ and $(P')$ is related to necessary conditions for the existence of positive solutions. It is easily seen that $a^+ \not \equiv 0$ is such a necessary condition for both problems. Furthermore, $\int_{\Omega}a<0$ is necessary for the existence of a positive solution of $(P')$, as shown in \cite[Lemma 2.1]{BPT}.
Under a stronger regularity condition on $a$, the existence of a nontrivial
nonnegative solution of $(P)$ has been proved in \cite{bandle,PT}. In
addition, the existence of a nontrivial nonnegative solution of $(P^{\prime})$
has been obtained in \cite{BPT} (see also \cite{alama}), under the additional
condition that $\int_{\Omega}a<0$. Furthermore, the authors in \cite{BPT} also
proved that the latter condition is necessary for the existence of positive
solutions for $(P^{\prime})$, if $f\in C^{1}(0,\infty)$ and $f^{\prime}(s)>0$
for $s>0$.
%\marginpar{\textit{22oct2016}}
%Via variational methods, it is easy to show the existence of a nontrivial non-negative weak solution of $(P)$, which turns to be a solution in the above sense. Still by variational methods, one can show that if then $(P')$ has a \textit{nontrivial non-negative }solution.
However, due to the non-Lipschitzian character of $f$ at $s=0$ and the change
of sign in $a$, neither the strong maximum principle nor Hopf's Lemma applies
to $(P)$ and $(P^{\prime})$. As a consequence, one can't deduce the positivity
of nontrivial nonnegative solutions of $(P)$ or $(P^{\prime})$.

In fact, the existence of a positive solution for these indefinite sublinear
problems is a delicate issue and very few papers in the literature have
addressed this question. Regarding $(P)$, when $f\left(  s\right)  =s^{q}$, it
was first proved in \cite{jesusultimo} that if the unique solution $\varphi$
of the linear problem%
\begin{equation}
\left\{
\begin{array}
[c]{lll}%
-\Delta\varphi=a(x) & \mathrm{in} & \Omega,\\
\varphi=0 & \mathrm{on} & \partial\Omega,
\end{array}
\right.  \label{jesuss}%
\end{equation}
is such that $\varphi>0$ in $\Omega$, then $(P)$ has a positive solution
(which may \textit{not} belong to the interior of the positive cone). This
condition, however, is not sharp, since one can find a function $a$ such that
$(P)$ possesses a positive solution but the corresponding $\varphi$ satisfies
that $\varphi<0$ in $\Omega$ (see \cite[Section 1]{nodea}). Later on, in the
aforementioned article \cite{nodea}, the authors considered the same problem
in the one-dimensional and radial cases, providing several sufficient
conditions on $a$ (as well as some necessary conditions) for the existence of
a positive solution of $(P)$. Some of these results were then extended to the
case of a general bounded domain in \cite{ans}. We point out that in all these
papers the only tool used was essentially the well known sub and
supersolutions method in the presence of weak and well-ordered sub and
supersolutions (see e.g. \cite{du}).
%For results on \textit{dead core} solutions of $(P)$ (with $f\left(  s\right)  =s^{q}$) we refer to \cite{bandle,PT}.

On the other hand, for the Neumann problem $(P^{\prime})$, even with $f\left(
s\right)  =s^{q}$, to the best of our knowledge, no \textbf{sufficient}
conditions for the existence of positive solutions are known.
%\marginpar{\textit{22oct2016}}

In this article, we shall not only prove that in some cases $(P^{\prime})$ and
$(P)$ admit positive solutions, but even more, that \textbf{every} nontrivial
nonnegative solution of $(P)$ and $(P^{\prime})$ is a positive solution. This
will be done using a continuity argument from \cite{J} (see also \cite{K}),
where the author proves the existence of a positive solution for the problem
\begin{equation}
\left\{
\begin{array}
[c]{lll}%
-\Delta u=u^{p}+f(x) & \mathrm{in} & \Omega,\\
u=0 & \mathrm{on} & \partial\Omega,
\end{array}
\right.  \label{jj}%
\end{equation}
where $1<p<\frac{N+2}{N-2}$, $N\geq3$, and $f\in L^{s}(\Omega)$, with
$s>\frac{N}{2}$. Under a smallness condition on $f$ (which may change sign),
the author shows that this problem has a mountain-pass solution $u_{f}$ which
depends continuously on $f$, in the sense that, up to a subsequence,
$u_{f}\rightarrow u_{0}$ in $\mathcal{C}^{1}(\overline{\Omega})$ as
$f\rightarrow0$ in $L^{s}(\Omega)$, where $u_{0}$ is a nontrivial nonnegative
solution of \eqref{jj} with $f\equiv0$. Furthermore, by the strong maximum
principle and Hopf's Lemma, $u_{0}$ lies in the interior of the positive cone
of $\mathcal{C}^{1}(\overline{\Omega})$, and consequently so does $u_{f}$ if
$f$ is close enough to zero. We shall exploit this idea, dealing now with a
class of sublinear problems and deducing the positivity of not only one
solution, but every nontrivial nonnegative solution.

Roughly speaking, we shall see that the positivity of nontrivial nonnegative
solutions can be recovered if $(P)$ or $(P^{\prime})$ are somehow sufficiently
close to a problem the strong maximum principle applies to. This situation
occurs, for instance, if the negative part of $a$ is small enough (for $(P)$)
or if $f(s)=s^{q}$ with $q$ close enough to $1$ (for both $(P)$ and
$(P^{\prime})$). We rely here on the fact that the strong maximum principle
applies to $(P)$ and $(P^{\prime})$ if either $a\geq0$ or $f(s)=s$.

We set $a^{\pm}:=\max(\pm a,0)$. Observe that the assumption that $a$ changes
sign means that $|\text{supp }a^{\pm}|>0$, where $|A|$ stands for the Lebesgue
measure of $A\subset\mathbb{R}^{N}$. We denote by $\Omega_{+}$ the largest
open subset of $\Omega$ where $a>0$ \textit{a.e.}, and assume that
\[
\Omega_{+}\text{ has finitely many connected components and }|(\text{supp
}a^{+})\setminus\Omega_{+}|=0.\leqno{(H_2)}
\]
In particular, we see that $\Omega_{+}$ is nonempty.

The above condition will be used to deduce that nontrivial nonnegative
solutions of $(P)$ and $(P^{\prime})$ are positive in a connected component of
$\Omega_{+}$, and consequently uniformly bounded away from zero therein (see
Lemma \ref{l1}). To this end, we shall also assume the following technical
condition, which is related to the use of the strong maximum principle:
\[
K_{s_{0}}:=\inf_{0\leq t<s\leq s_{0}}\frac{f(s)-f(t)}{s-t}>-\infty
,\quad\text{for all }s_{0}>0.\leqno{(H_3)}
\]
Note in particular that this condition is satisfied, for instance, if $f$ is
nondecreasing (in which case $K_{s_{0}}\geq0$), and in particular,
$f(s)=s^{q}$ with $0\leq q<1$. 

Our positivity results for $(P)$ shall provide us with solutions that lie in
the interior of the positive cone of $\mathcal{C}_{0}^{1}(\overline{\Omega
}):=\{u\in\mathcal{C}^{1}(\overline{\Omega}):u=0\text{ on }\partial\Omega\}$,
which is denoted by
\[
\mathcal{P}_{D}^{\circ}:=\left\{  u\in\mathcal{C}_{0}^{1}(\overline{\Omega
}):u>0\text{ in }\Omega\text{, \ and }\frac{\partial u}{\partial\nu}<0\text{
on }\partial\Omega\right\}  .
\]
Regarding $(P^{\prime})$, we shall obtain solutions that belong to
\[
\mathcal{P}_{N}^{\circ}:=\left\{  u\in\mathcal{C}^{1}(\overline{\Omega
}):u>0\text{ on }\overline{\Omega}\right\}  .
\]
Note that a positive solution of $(P)$ (respect. $(P^{\prime})$) \textit{need
not} belong to $\mathcal{P}_{D}^{\circ}$ (respect. $\mathcal{P}_{N}^{\circ}$),
as shown in Proposition \ref{pex} below.
%\marginpar{\textit{22oct2016}}

We state now our main results.

\begin{theorem}
\label{t1} Assume $(H_{1})$, $(H_{2})$, and $(H_{3})$. Then there exists
$\delta>0$ (possibly depending on $a^{+}$) such that every nontrivial
nonnegative solution of $(P)$ belongs to $\mathcal{P}_{D}^{\circ}$ if $\Vert
a^{-}\Vert_{L^{r}\left(  \Omega\right)  }<\delta$.
\end{theorem}

\begin{remark}
\label{rem:nec} As already mentioned, if $f\in C^{1}(0,\infty)$ and
$f^{\prime}>0$ in $(0,\infty)$, then the condition $\int_{\Omega}a<0$ is
necessary for the existence of positive solutions of $(P^{\prime})$, cf.
\cite[Lemma 2.1]{BPT}. In view of this fact, we can't expect an analogue of
Theorem \ref{t1} for $(P^{\prime})$.
\end{remark}

In the case that $f$ is a power, we write $(P)$ and $(P^{\prime})$ as
\[
\left\{
\begin{array}
[c]{lll}%
-\Delta u=a(x)u^{q} & \mathrm{in} & \Omega,\\
u=0 & \mathrm{on} & \partial\Omega,
\end{array}
\right.  \leqno{(P_D)}
\]
and
\[
\left\{
\begin{array}
[c]{lll}%
-\Delta u=a(x)u^{q} & \mathrm{in} & \Omega,\\
\frac{\partial u}{\partial\nu}=0 & \mathrm{on} & \partial\Omega.
\end{array}
\right.  \leqno{(P_N)}
\]

\begin{theorem}
\label{t2} Assume $(H_{2})$. Then, given $a\in L^{r}\left(  \Omega\right)  $,
there exists $q_{0}\in(0,1)$ such that every nontrivial nonnegative solution
of $(P_{D})$ belongs to $\mathcal{P}_{D}^{\circ}$ if $q_{0}<q<1$.
\end{theorem}

As a consequence of Theorems \ref{t1} and \ref{t2}, we derive the following
existence and uniqueness results:

\begin{corollary}
\label{c1} Under the conditions of Theorem \ref{t1}, let $\Vert a^{-}%
\Vert_{L^{r}\left(  \Omega\right)  }<\delta$. Assume in addition that
$f\in\mathcal{C}^{1}(0,\infty)$, $f^{\prime}$ is nonincreasing in $(0,\infty)$
and $\int_{0}^{t}\frac{1}{f(s)}ds<\infty$ for $t>0$. Then $(P)$ has a solution
in $\mathcal{P}_{D}^{\circ}$ and has no other nontrivial nonnegative solutions.
\end{corollary}

\begin{corollary}
\label{c2} Under the assumptions of Theorem \ref{t2}, let $q_{0}<q<1$. Then
$(P_{D})$ has a solution in $\mathcal{P}_{D}^{\circ}$ and has no other
nontrivial nonnegative solutions.
\end{corollary}

\begin{remark}
Let us mention that if $f$ is nondecreasing and $k_{1}s^{q}\leq f\left(
s\right)  \leq k_{2}s^{q}$ for some $k_{1},k_{2}>0$ and all $s\geq0$, as a
consequence of \cite[Theorem 3.1]{ans}, one can deduce the existence of a
solution $u\in\mathcal{P}_{D}^{\circ}$ for $(P)$ if $a^{-}$ is sufficiently
small with respect to $a^{+}$. On the other side, in the one-dimensional and
radial cases one can derive the existence of a positive solution of $(P_{D})$
(but not necessarily belonging to $\mathcal{P}_{D}^{\circ}$) provided that $q$
is close enough to $1$ (cf. \cite{nodea}). In this sense, Corollaries \ref{c1}
and \ref{c2} are consistent with the existence results from \cite{nodea,ans}.
\end{remark}

For the Neumann problem $(P_{N})$, we establish the following analogue of
Theorem \ref{t2}: \medskip

\begin{theorem}
\label{t3} Assume $(H_{2})$. Then, given $a\in L^{r}\left(  \Omega\right)  $,
there exists $q_{0}\in(0,1)$ such that every nontrivial nonnegative solution
of $(P_{N})$ belongs to $\mathcal{P}_{N}^{\circ}$ if $q_{0}<q<1$.
\end{theorem}

\begin{corollary}
\label{c3} Under the assumptions of Theorem \ref{t3}, let $q_{0}<q<1$. Assume
in addition that $\int_{\Omega}a<0$. Then $(P_{N})$ has a solution in
$\mathcal{P}_{N}^{\circ}$ and has no other nontrivial nonnegative solutions.
\end{corollary}

Our next result concerns the sets
\begin{align*}
\mathcal{A}_{D}  &  :=\{q\in(0,1):\text{any nontrivial nonnegative solution of
}(P_{D})\text{ lies in }\mathcal{P}_{D}^{\circ}\},\\
\mathcal{A}_{N}  &  :=\{q\in(0,1):\text{any nontrivial nonnegative solution of
}(P_{N})\text{ lies in }\mathcal{P}_{N}^{\circ}\}.
\end{align*}

\begin{theorem}
\label{ta} Assume $(H_{2})$ and fix $a\in L^{r}\left(  \Omega\right)  $. Then:

\begin{enumerate}
\item[(i)] $\mathcal{A}_{D}$ is a nonempty open interval.

\item[(ii)] $\mathcal{A}_{N}$ is a nonempty open interval under the condition
$\int_{\Omega}a<0$.
\end{enumerate}

In particular, there exists $q_{0} \in[0,1)$ such that $\mathcal{A}_{D}%
=(q_{0},1)$, and a similar characterization holds for $\mathcal{A}_{N}$.
\end{theorem}

\begin{remark}
Although Theorem \ref{ta} states that (under $(H_{2})$) the sets
$\mathcal{A}_{D}$ and $\mathcal{A}_{N}$ (assuming $\int_{\Omega}a<0$) are
always nonempty, as a consequence of Proposition \ref{prop} below, we shall
see that given any $q\in(0,1)$, we may find $a$ in such a way that $(P_{D})$
(respect. $(P_{N})$) has nontrivial nonnegative solutions that do not belong
to $\mathcal{P}_{D}^{\circ}$ (respect. $\mathcal{P}_{N}^{\circ}$). This fact
shows that $\mathcal{A}_{D}$ and $\mathcal{A}_{N}$ can be made arbitrarily
small by choosing $a$ in a suitable way.
%Let us also mention
%that the aforementioned proposition shows that, given any $q\in\left(
%0,1\right)$, one can find a weight $a$ such that $(P_{D})$ has nontrivial non-negative \textit{dead core
%} solutions (i.e., nontrivial solutions vanishing somewhere in $\Omega$), even if such
%$q$ is close to $1$ (see Remark \ref{core} below).

\end{remark}

The rest of the paper is organized as follows. In the next section we prove
some auxiliary results concerning $(P)$ and $(P^{\prime})$, whereas in Section
3 we supply the proofs of our main results. Finally, in Section 4 we apply
some of our theorems to derive positivity results (as well as existence and
multiplicity results of positive solutions) for indefinite concave-convex type
problems. \medskip

\section{Preliminary results}

\medskip Let us fix the notation to be used in the sequel.

Given $m\in L^{r}\left(  \Omega\right)  $, $r>N$, and an open set
$B\subseteq\Omega$ such that $m^{+}\not \equiv 0$ in $B$, we denote by
$\lambda_{1}(m,B)$ the first positive eigenvalue of the problem
\[
\left\{
\begin{array}
[c]{lll}%
-\Delta\phi=\lambda m\left(  x\right)  \phi & \mathrm{in} & B,\\
\phi=0 & \mathrm{on} & \partial B.
\end{array}
\right.
\]
We shall deal with several norms, which will be denoted as follows: $\Vert
u\Vert_{L^{r}(\Omega)}:=(\int_{\Omega}|u|^{r})^{\frac{1}{r}}$, $\Vert
u\Vert_{H_{0}^{1}(\Omega)}:=(\int_{\Omega}|\nabla u|^{2})^{\frac{1}{2}}$ and
$\Vert u\Vert_{H^{1}(\Omega)}:=(\int_{\Omega}\left(  |\nabla u|^{2}%
+u^{2}\right)  )^{\frac{1}{2}}$.\smallskip

To begin with, we provide several useful lower bounds for nontrivial
nonnegative solutions of $(P)$ and $(P^{\prime})$.
%Its proof is inspired by the proof of
%\cite[Proposition 3.5]{AM}.

\begin{lemma}
\label{l0} Assume $(H_{2})$ and let $u$ be a nontrivial nonnegative solution
of $(P)$ or $(P^{\prime})$. Then there exists a subdomain $\Omega^{\prime
}\subset\Omega_{+}$ such that $u>0$ in $\Omega^{\prime}$.
\end{lemma}

\textbf{Proof}. If $u$ is a nontrivial nonnegative solution of $(P)$ then it
satisfies
\[
0<\int_{\Omega}|\nabla u|^{2}=\int_{\Omega}a(x)f(u)u\leq\int_{\text{supp
}a^{+}}a^{+}(x)f(u)u=\int_{\Omega_{+}}a^{+}(x)f(u)u,
\]
where we used the assumption that $|(\text{supp }a^{+})\setminus\Omega_{+}%
|=0$. It follows that $u\not \equiv 0$ in $\Omega_{+}$, and consequently $u>0$
in some subdomain of $\Omega_{+}$. The same argument applies if $u$ is a
nontrivial nonnegative solution of $(P^{\prime})$, since $u$ can't be a
constant. $\blacksquare$\newline

\begin{lemma}
\label{l1} Assume that $(H_{3})$ holds, $\displaystyle\lim_{s\rightarrow0^{+}%
}s^{-1}f(s)=\infty$ and $\Omega^{\prime}\not =\emptyset$ is a subdomain of
$\Omega_{+}$. Then, for any open ball $B$ such that $\overline{B}\subset
\Omega^{\prime}$ there exists a function $\psi\in W^{2,r}\left(  B\right)  $
such that $u\geq\psi>0$ in $B$ for every nontrivial nonnegative supersolution
of
\begin{equation}
-\Delta u=a(x)f(u)\text{\quad in }\Omega^{\prime}. \label{ss}%
\end{equation}

\end{lemma}

\textbf{Proof}. Let $u$ be a nontrivial nonnegative supersolution of
\eqref{ss} and $B$ be an open ball such that $\overline{B}\subset
\Omega^{\prime}$. Then $a\geq0$ and $a\not \equiv 0$ in $B$. Let $\phi\in
W^{2,r}\left(  B\right)  \cap W_{0}^{1,r}\left(  B\right)  $ be a positive
eigenfunction associated to $\lambda_{1}(a,B)$, with $\left\Vert
\phi\right\Vert _{\infty}=1$. We observe that for all $\varepsilon>0$
sufficiently small it holds that
\[
-\Delta\left(  \varepsilon\phi\right)  \leq a(x)f(\varepsilon\phi)\text{\quad
in }B.
\]
Indeed, note that $-\Delta(\varepsilon\phi)=\varepsilon\lambda_{1}%
(a,B)a(x)\phi$ in $B$. Hence, it is enough to check that $\varepsilon
\lambda_{1}(a,B)a(x)\phi\leq a(x)f(\varepsilon\phi)$, i.e.
\[
\lambda_{1}(a,B)\leq\frac{f(\varepsilon\phi)}{\varepsilon\phi}\text{\quad in
}B.
\]
Since $\frac{f(s)}{s}\rightarrow\infty$ as $s\rightarrow0^{+}$ and $\phi$ is
bounded, we see that there exists $\varepsilon_{0}>0$ such that the above
inequality holds for all $0<\varepsilon\leq\varepsilon_{0}$.

To conclude the proof we show that $u\geq\varepsilon_{0}\phi$ in $B$ (note
that $\varepsilon_{0}$ does not depend on $u$). Indeed, suppose this is not
true. Since $\Omega^{\prime}$ is connected, it follows from the strong maximum
principle that $u>0$ in $\Omega^{\prime}$, so that $u>0$ in $\overline{B}$.
Moreover, $\phi=0$ on $\partial B$, so there exists $s\in\left(  0,1\right)  $
such that $u\geq s\varepsilon_{0}\phi$ in $B$ and $u\left(  x_{0}\right)
=s\varepsilon_{0}\phi\left(  x_{0}\right)  $ for some $x_{0}\in B$. Setting
$s_{0}:=\left\Vert u\right\Vert _{L^{\infty}\left(  \Omega\right)  }$ and
$M(x):=\left\vert K_{s_{0}}\right\vert a(x)$, where $K_{s_{0}}$ is given by
$(H_{3})$, one can see that the map $s\rightarrow M(x)s+a\left(  x\right)
f\left(  s\right)  $ is nondecreasing for all $s\in\left(  0,\left\Vert
u\right\Vert _{\infty}\right)  $ and $a.e.$ $x\in B$. Then,%
\begin{align*}
&  -\Delta(u-s\varepsilon_{0}\phi)+M\left(  x\right)  \left(  u-s\varepsilon
_{0}\phi\right) \\
&  \geq M\left(  x\right)  \left(  u-s\varepsilon_{0}\phi\right)  +a(x)\left(
f(u)-f(s\varepsilon_{0}\phi)\right)  \geq0\text{\quad in }B,
\end{align*}
and $u>s\varepsilon_{0}\phi$ on $\partial B$. Therefore, the strong maximum
principle (e.g. \cite{tru}) says that $u>s\varepsilon_{0}\phi$ in $B$, which
is a contradiction. The proof is complete. $\blacksquare$

\begin{remark}
\label{rl1} Let us point out that the proof of Lemma \ref{l1} is the only
instance where $(H_{3})$ is employed.
\end{remark}

We prove now that $\left\Vert u\right\Vert _{H_{0}^{1}(\Omega)}\geq C$ for any
nontrivial nonnegative solution of $(P)$, for some constant $C>0$ independent
of $a^{-}$.

\begin{lemma}
\label{l2} Assume that $(H_{2})$ and $(H_{3})$ hold, and $\displaystyle\lim
_{s\rightarrow0^{+}}s^{-1}f(s)=\infty$. Then there exists a constant $C>0$
such that $\Vert u\Vert_{H_{0}^{1}(\Omega)}\geq C$ for every nontrivial
nonnegative solution $u$ of $(P)$. Moreover, $C$ does not depend on $a^{-}$.
\end{lemma}

\textbf{Proof}. Assume by contradiction that there exists a sequence $\left\{
u_{n}\right\}  $ of solutions of $(P)$ with $u_{n}\rightarrow0$ in $H_{0}%
^{1}(\Omega)$. Then $u_{n}\rightarrow0$ in $L^{2}(\Omega)$, and, up to a
subsequence, we have $u_{n}\rightarrow0$ \textit{a.e.} in $\Omega$. By Lemma
\ref{l0}, we know that any nontrivial nonnegative solution $u$ of $(P)$ is
positive in some subdomain of $\Omega_{+}$. Thus, since $\Omega_{+}$ has
finitely many connected components, we may assume that, for all $n\in
\mathbb{N}$, $u_{n}>0$ in some fixed subdomain $\widetilde{\Omega}%
\subset\Omega_{+}$. However, by Lemma \ref{l1}, we have $u_{n}\geq\psi>0$ in
some open ball $B\subset\widetilde{\Omega}$, so we reach a contradiction.
$\blacksquare$\newline

Next we get an \textit{a priori} bound from below for nontrivial nonnegative
solutions of either $(P_{D})$ or $(P_{N})$. We remark that this estimate does
not depend on $q$.

\begin{lemma}
\label{l12} Assume that $\Omega^{\prime}\not =\emptyset$ is a subdomain of
$\Omega_{+}$ such that $\lambda_{1}(a,\Omega^{\prime})<1$. Then there exists
an open set $B$ such that $\overline{B}\subset\Omega^{\prime}$ and a function
$\phi\in W^{2,r}(B)$ such that $u\geq\phi>0$ in $B$, for every nontrivial
nonnegative supersolution of
\begin{equation}
-\Delta u=a(x)u^{q}\text{\quad in }\Omega^{\prime}, \label{epr}%
\end{equation}
and for every $q\in(0,1)$.
\end{lemma}

\textbf{Proof}. Let $\Omega_{\delta}^{\prime}:=\{x\in\Omega^{\prime
}:dist(x,\partial\Omega^{\prime})>\delta\}$ for $\delta>0$. From the
variational characterization of $\lambda_{1}(a,B)$, we know that $\lambda
_{1}(a,\Omega_{\delta}^{\prime})\rightarrow\lambda_{1}(a,\Omega^{\prime})$ as
$\delta\rightarrow0^{+}$ (see e.g. Lemma 2.5 in \cite{DFG}). We fix
$\delta_{0}>0$ such that $\lambda_{1}(a,\Omega_{\delta_{0}}^{\prime})<1$ and
set $B:=\Omega_{\delta_{0}}^{\prime}$. Let $\phi>0$ with $\Vert\phi
\Vert_{\infty}=1$ be as in Lemma \ref{l1}, i.e. a solution of
\[
\left\{
\begin{array}
[c]{lll}%
-\Delta\phi=\lambda_{1}\left(  a,B\right)  a\left(  x\right)  \phi &
\mathrm{in} & B,\\
\phi=0 & \mathrm{on} & \partial B,
\end{array}
\right.
\]
and let
\[
0<\varepsilon\leq\varepsilon_{q}:=\frac{1}{\lambda_{1}\left(  a,B\right)
^{1/(1-q)}}.
\]
Observe that since $\lambda_{1}\left(  a,B\right)  <1$ we have $\varepsilon
_{q}\geq1$ for all $q$. Then, taking into account that $0<\phi\leq1$ and the
definition of $\varepsilon_{q}$, we derive that
\[
-\Delta\left(  \varepsilon\phi\right)  =\lambda_{1}\left(  a,B\right)
a\left(  x\right)  \varepsilon\phi\leq a(x)\left(  \varepsilon\phi\right)
^{q}\text{\quad in }B
\]
for all $0<\varepsilon\leq\varepsilon_{q}$.

Now, given any nontrivial nonnegative supersolution $u$ of (\ref{epr}), we can
argue as in the last paragraph of the proof of Lemma \ref{l1} to infer that
$u\geq\phi$ in $B$ for any $q\in(0,1)$. This concludes the proof of the lemma.
$\blacksquare$ \newline

The proof of the next estimates are similar to the one of Lemma \ref{l2} (we
use now Lemma \ref{l12}), so we omit it.

\begin{lemma}
\label{l3} Assume $(H_{2})$ and $0<q<1$. Then there exists a constant $C>0$
such that:

\begin{enumerate}
\item[(i)] $\Vert u\Vert_{H_{0}^{1}(\Omega)}\geq C$ for every nontrivial
nonnegative solution $u$ of $(P_{D})$.

\item[(ii)] $\Vert u\Vert_{H^{1}(\Omega)}\geq C$ for every nontrivial
nonnegative solution $u$ of $(P_{N})$.
\end{enumerate}

Moreover, $C$ does not depend on $q$.
\end{lemma}

To end this section, we prove some results on the sets $\mathcal{A}_{D}$ and
$\mathcal{A}_{N}$.

\begin{lemma}
\label{l6} Assume $(H_{2})$. If $q_{0}\in\mathcal{A}_{D}$, then $\left(
q_{0},\frac{1}{2-q_{0}}\right)  \subset\mathcal{A}_{D}$.
\end{lemma}

\textbf{Proof}. Assume to the contrary that $q_{0}\in\mathcal{A}_{D}$,
$q_{0}<q<\frac{1}{2-q_{0}}$, but $q\not \in \mathcal{A}_{D}$. It follows that
there exists a nontrivial nonnegative solution $u$ of $(P_{D})$ such that
$u\not \in \mathcal{P}_{D}^{\circ}$. Let
\[
\beta:=\frac{1-q}{1-q_{0}}\in\left(  0,1\right)  ,\qquad\gamma:=\frac{q-q_{0}%
}{1-q}>0,
\]
and consider the auxiliary problem
\begin{equation}%
\begin{cases}
-\Delta w=\beta a(x)w^{-\gamma}(w^{1/\beta}-\varepsilon)^{q} &
\mbox{in $\Omega$},\\
w\geq\varepsilon^{\beta} & \mbox{in $\Omega$},\\
w=\varepsilon^{\beta} & \mbox{on $\partial \Omega$},
\end{cases}
\label{pa}%
\end{equation}
with $0<\varepsilon\leq1$. Equivalently, putting $\hat{w}:=w-\varepsilon
^{\beta}$, we consider
\begin{equation}%
\begin{cases}
-\Delta\hat{w}=\beta a(x)(\hat{w}+\varepsilon^{\beta})^{-\gamma}\{(\hat
{w}+\varepsilon^{\beta})^{1/\beta}-\varepsilon\}^{q} & \mbox{in $\Omega$},\\
w\geq0 & \mbox{in $\Omega$},\\
w=0 & \mbox{on $\partial \Omega$}.
\end{cases}
\label{pa2}%
\end{equation}
The limiting problem as $\varepsilon\rightarrow0^{+}$ is understood as
\begin{equation}%
\begin{cases}
-\Delta w=\beta a(x)w^{q_{0}} & \mbox{in $\Omega$},\\
w\geq0 & \mbox{in $\Omega$},\\
w=0 & \mbox{on $\partial \Omega$}.
\end{cases}
\label{lp}%
\end{equation}
Since $q_{0}\in\mathcal{A}_{D}$, any nontrivial nonnegative solution of
\eqref{lp} belongs to $\mathcal{P}_{D}^{\circ}$.

In the sequel, we shall obtain a solution $w_{\varepsilon}$ of \eqref{pa} and
show that, as $\varepsilon\rightarrow0$, $w_{\varepsilon}$ converges (up to a
subsequence) to a nontrivial nonnegative solution of \eqref{lp} that does not
belong to $\mathcal{P}_{D}^{\circ}$. This will provide us with a
contradiction. We divide the rest of the proof in several steps:\newline

\textbf{Step 1:} Construction of a weak supersolution of \eqref{pa}.\newline

We note that $\psi=\psi_{\varepsilon}:=(u+\varepsilon)^{\beta}$ is a
supersolution of \eqref{pa}. Indeed, since $1-\beta=\gamma\beta$, by direct
computations we have that
\[
-\Delta\psi\geq\beta a(x)\psi^{-\gamma}(\psi^{1/\beta}-\varepsilon)^{q}%
\quad\mbox{in $\Omega$},
\]
and $\psi=\varepsilon^{\beta}$ on $\partial\Omega$, as desired. \newline

\textbf{Step 2:} Construction of a weak subsolution of \eqref{pa}.\newline

%Next, we shall construct a subsolution of \eqref{pa}.
By Lemma \ref{l0}, there exists a ball $B$ such that $a\geq0$, $a\not \equiv
0$ $a.e.$ in $B$ and $u>0$ in $\overline{B}$. Let $\phi\in W^{2,r}\left(
B\right)  \cap W_{0}^{1,r}\left(  B\right)  $ be a positive eigenfunction
associated to $\lambda_{1}(a,B)$, with $\left\Vert \phi\right\Vert _{\infty
}=1$, and extend $\phi$ to $\overline{\Omega}$ by $\phi=0$ in $\overline
{\Omega}\setminus\overline{B}$. Given $0<\delta\leq1$, we set%
\begin{equation}
\varphi_{\delta,\varepsilon}:=\left\{
\begin{array}
[c]{ll}%
\delta\phi+\varepsilon^{\beta} & \text{in }B,\\
\varepsilon^{\beta} & \text{in }\overline{\Omega}\setminus\overline{B}.
\end{array}
\right.  \label{zzz}%
\end{equation}
We observe that
\begin{align*}
&  -\Delta\varphi_{\delta,\varepsilon}-\beta a(x)\varphi_{\delta,\varepsilon
}^{-\gamma}(\varphi_{\delta,\varepsilon}^{1/\beta}-\varepsilon)^{q}\\
&  \leq a\left(  x\right)  \left(  \lambda_{1}\left(  a,B\right)  \delta
\phi-\beta(\delta\phi+\varepsilon^{\beta})^{-\gamma}\{(\delta\phi
+\varepsilon^{\beta})^{1/\beta}-\varepsilon\}^{q}\right)  \quad\mbox{in $B$}.
\end{align*}
We claim that there exists $c_{0}>0$, independent of $x\in B$, such that for
$\varepsilon\in(0,1]$ and $\delta\in(0,1]$ we have
\[
(\delta\phi+\varepsilon^{\beta})^{-\gamma}\{(\delta\phi+\varepsilon^{\beta
})^{1/\beta}-\varepsilon\}^{q}\geq c_{0}(\delta\phi)^{q_{0}+\gamma}%
\quad\mbox{in $B$}.
\]
Indeed, since $q/\beta=q_{0}+\gamma$, we note that for $x\in B$,
\begin{align*}
\frac{(\delta\phi+\varepsilon^{\beta})^{-\gamma}\{(\delta\phi+\varepsilon
^{\beta})^{1/\beta}-\varepsilon\}^{q}}{(\delta\phi)^{q_{0}+\gamma}}  &
\geq\frac{(\delta\phi)^{q/\beta}}{(\delta\phi+\varepsilon^{\beta})^{\gamma
}(\delta\phi)^{q_{0}+\gamma}}\\
&  =\frac{1}{(\delta\phi+\varepsilon^{\beta})^{\gamma}}\geq\frac{1}{2^{\gamma
}}=:c_{0},
\end{align*}
as desired. Here, we have used the fact that if $\alpha>1$ then $(s+t)^{\alpha
}\geq s^{\alpha}+t^{\alpha}$ for $t,s\geq0$. Thus the claim is proved. It
follows that for $\delta>0$ small enough, we have that
\begin{align*}
-\Delta\varphi_{\delta,\varepsilon}-\beta a(x)\varphi_{\delta,\varepsilon
}^{-\gamma}\left(  \varphi_{\delta,\varepsilon}^{1/\beta}-\varepsilon\right)
^{q}  &  \leq a\left(  x\right)  \left(  \lambda_{1}\left(  a,B\right)
\delta\phi-c_{0}\beta(\delta\phi)^{q_{0}+\gamma}\right) \\
&  \leq0\quad\mbox{in $B$},
\end{align*}
since the assumption $q<\frac{1}{2-q_{0}}$ implies $q_{0}+\gamma<1$. Note that
$\delta$ is determined uniformly in $\varepsilon\in(0,1]$. Thus, employing the
divergence theorem as stated e.g. in \cite{cuesta}, p. 742, we deduce that
$\varphi_{\delta,\varepsilon}$ is a weak subsolution of \eqref{pa}.\newline

\textbf{Step 3:} The subsolution and the supersolution of \eqref{pa} are
well-ordered.\newline

We shall see that, choosing $\delta$ and $\varepsilon$ adequately,
$(u+\varepsilon)^{\beta}$ and $\varphi_{\delta,\varepsilon}$ are well-ordered,
i.e., $(u+\varepsilon)^{\beta}\geq\varphi_{\delta,\varepsilon}$ in $\Omega$.
We assert that there exist $\delta_{1},\varepsilon_{1}>0$ such that if
$\varepsilon\in(0,\varepsilon_{1})$, then
\[
\left(  u+\varepsilon\right)  ^{\beta}\geq\delta_{1}\phi+\varepsilon^{\beta
}\quad\mbox{in $B$}.
\]
Indeed, if we fix $\varepsilon_{1}$ and $\delta_{1}$ such that%
\[
\varepsilon_{1}^{\beta}\leq\frac{1}{2}(\min_{\overline{B}}u)^{\beta}%
,\qquad\delta_{1}\leq\frac{1}{4}(\min_{\overline{B}}u)^{\beta},
\]
then it is clear that
\begin{align*}
(u+\varepsilon)^{\beta}  &  \geq\frac{1}{2}u^{\beta}+\frac{1}{2}%
\varepsilon^{\beta}\geq\frac{1}{2}\{(\min_{\overline{B}}u)^{\beta}%
-\varepsilon^{\beta}\}+\varepsilon^{\beta}\\
&  \geq\frac{1}{4}(\min_{\overline{B}}u)^{\beta}+\varepsilon^{\beta}\geq
\delta_{1}\phi+\varepsilon^{\beta}\quad\mbox{in $B$}.
\end{align*}

Hence, for every $\varepsilon\in(0,\varepsilon_{1})$, the method of weak sub
and supersolutions (see e.g. \cite[Theorem 4.9]{du}) gives us some
$w_{\varepsilon}\in H_{0}^{1}\left(  \Omega\right)  \cap L^{\infty}\left(
\Omega\right)  $ solution of \eqref{pa}, with
\begin{equation}
\varphi_{\delta_{1},\varepsilon}\leq w_{\varepsilon}\leq(u+\varepsilon
)^{\beta}\quad\mbox{in $\Omega$}. \label{ssmm}%
\end{equation}
Furthermore, by standard regularity arguments, $w_{\varepsilon}\in
W^{2,r}\left(  \Omega\right)  \cap W_{0}^{1,r}\left(  \Omega\right)
$.\newline

\textbf{Step 4:} The limiting behavior of $w_{\varepsilon}$ as $\varepsilon
\rightarrow0^{+}$.\newline

We convert $w_{\varepsilon}$ to \eqref{pa2} by $\hat{w_{\varepsilon}%
}=w_{\varepsilon}-\varepsilon^{\beta}$, so that $\hat{w}_{\varepsilon}=0$ on
$\partial\Omega$. Thus, we deduce that
\begin{align*}
\int_{\Omega}|\nabla\hat{w}_{\varepsilon}|^{2}  &  =\beta\int_{\Omega}%
a(x)\hat{w}_{\varepsilon}(\hat{w}_{\varepsilon}+\varepsilon^{\beta})^{-\gamma
}((\hat{w}_{\varepsilon}+\varepsilon^{\beta})^{1/\beta}-\varepsilon)^{q}\\
&  \leq C\int_{\Omega}\hat{w}_{\varepsilon}(\hat{w}_{\varepsilon}%
+\varepsilon^{\beta})^{-\gamma}(\hat{w}_{\varepsilon}+\varepsilon^{\beta
})^{q/\beta}\\
&  =C\int_{\Omega}\hat{w}_{\varepsilon}(\hat{w}_{\varepsilon}+\varepsilon
^{\beta})^{q_{0}}.
\end{align*}
Since we see from \eqref{ssmm} that $\Vert\hat{w}_{\varepsilon}\Vert
_{L^{\infty}(\Omega)}\leq C$ as $\varepsilon\rightarrow0^{+}$, we infer that
$\Vert\hat{w}_{\varepsilon}\Vert_{H_{0}^{1}(\Omega)}$ is bounded as
$\varepsilon\rightarrow0^{+}$. It follows that, up to a subsequence, $\hat
{w}_{\varepsilon}\rightharpoonup\hat{w}_{0}$ in $H_{0}^{1}(\Omega)$, and
$\hat{w}_{\varepsilon}\rightarrow\hat{w}_{0}$ $a.e.$ in $\Omega$ for some
$\hat{w}_{0}\in H_{0}^{1}(\Omega)$. Also, since $\hat{w}_{\varepsilon}$ is a
weak solution of \eqref{pa2}, we note that
\[
\int_{\Omega}\nabla\hat{w}_{\varepsilon}\nabla v=\beta\int_{\Omega}%
a(x)(\hat{w}_{\varepsilon}+\varepsilon^{\beta})^{-\gamma}\{\left(  \hat
{w}_{\varepsilon}+\varepsilon^{\beta}\right)  ^{1/\beta}-\varepsilon
\}^{q}v,\quad\forall v\in\mathcal{C}_{0}^{1}\left(  \Omega\right)  .
\]
So, from the fact that $\hat{w}_{\varepsilon}\rightharpoonup\hat{w}_{0}$ in
$H_{0}^{1}(\Omega)$, we get that
\[
\int_{\Omega}\nabla\hat{w}_{\varepsilon}\nabla v\rightarrow\int_{\Omega}%
\nabla\hat{w}_{0}\nabla v,\quad\forall v\in\mathcal{C}_{0}^{1}\left(
\Omega\right)  .
\]

On the other hand, recalling \eqref{ssmm} and that $-\gamma+q/\beta=q_{0}$, we
see that%

\begin{align*}
\left\vert a(x)(\hat{w}_{\varepsilon}+\varepsilon^{\beta})^{-\gamma}\{\left(
\hat{w}_{\varepsilon}+\varepsilon^{\beta}\right)  ^{1/\beta}-\varepsilon
\}^{q}v\right\vert  &  \leq C\left\vert a\left(  x\right)  \right\vert
(\hat{w}_{\varepsilon}+\varepsilon^{\beta})^{-\gamma}(\hat{w}_{\varepsilon
}+\varepsilon^{\beta})^{q/\beta}\\
&  =C\left\vert a\left(  x\right)  \right\vert (\hat{w}_{\varepsilon
}+\varepsilon^{\beta})^{q_{0}}\\
&  \leq C^{\prime}\left\vert a\left(  x\right)  \right\vert \in L^{r}\left(
\Omega\right)  .
\end{align*}
Therefore, the Lebesgue convergence theorem yields that%
\[
\int_{\Omega}a(x)(\hat{w}_{\varepsilon}+\varepsilon^{\beta})^{-\gamma
}\{\left(  \hat{w}_{\varepsilon}+\varepsilon^{\beta}\right)  ^{1/\beta
}-\varepsilon\}^{q}v\rightarrow\int_{\Omega}a(x)\hat{w}_{0}^{q_{0}}%
v,\quad\forall v\in\mathcal{C}_{0}^{1}\left(  \Omega\right)  .
\]
Indeed, if $w_{0}>0$, then
\[
(\hat{w}_{\varepsilon}+\varepsilon^{\beta})^{-\gamma}\{\left(  \hat
{w}_{\varepsilon}+\varepsilon^{\beta}\right)  ^{1/\beta}-\varepsilon
\}^{q}=w_{\varepsilon}^{-\gamma}(w_{\varepsilon}^{1/\beta}-\varepsilon
)^{q}\rightarrow w_{0}^{-\gamma}(w_{0}^{1/\beta})^{q}=w_{0}^{q_{0}},
\]
whereas if $w_{0}=0$, then%
\begin{align*}
(\hat{w}_{\varepsilon}+\varepsilon^{\beta})^{-\gamma}\{\left(  \hat
{w}_{\varepsilon}+\varepsilon^{\beta}\right)  ^{1/\beta}-\varepsilon\}^{q}  &
=w_{\varepsilon}^{-\gamma}(w_{\varepsilon}^{1/\beta}-\varepsilon)^{q}\\
&  \leq w_{\varepsilon}^{-\gamma}w_{\varepsilon}^{q/\beta}=w_{\varepsilon
}^{q_{0}}\rightarrow0.
\end{align*}

Summing up, we have obtained that
\[
\int_{\Omega}\nabla\hat{w}_{0}\nabla v=\beta\int_{\Omega}a(x)\hat{w}%
_{0}^{q_{0}}v=0,\quad\forall v\in\mathcal{C}_{0}^{1}(\Omega).
\]
This implies that $\hat{w}_{0}$ is a weak solution of \eqref{lp}. Now, from
\eqref{ssmm}, we recall that for $\varepsilon\in(0,\varepsilon_{1})$,
\[
\varphi_{\delta_{1},\varepsilon}-\varepsilon^{\beta}\leq\hat{w}_{\varepsilon
}\leq(u+\varepsilon)^{\beta}-\varepsilon^{\beta}\quad\mbox{in $\Omega$}.
\]
Therefore, passing to the limit as $\varepsilon\rightarrow0^{+}$, this
inequality provides
\[
\varphi_{\delta_{1},0}\leq\hat{w}_{0}\leq u^{\beta}\quad\mbox{in $\Omega$}.
\]
This means that $\hat{w}_{0}$ is a nontrivial nonnegative solution of
\eqref{lp}, but $\hat{w}_{0}\not \in \mathcal{P}_{D}^{\circ}$, since
$u\not \in \mathcal{P}_{D}^{\circ}$ by assumption. Hence we reach a
contradiction, and the proof is complete. $\blacksquare$

\begin{remark}
Lemma \ref{l6} also holds for $\mathcal{A}_{N}$. Indeed, we can prove Lemma
\ref{l6} for $\mathcal{A}_{N}$ with some minor modifications in the proof.
Assume $q_{0}\in\mathcal{A}_{N}$, $q_{0}<q<\frac{1}{2-q_{0}}$, but
$q\not \in \mathcal{A}_{N}$. It follows that there exists a nontrivial
nonnegative solution $u$ of $(P_{N})$ such that $u$ does not belong to
$\mathcal{P}_{N}^{\circ}$. The rest of the proof for $\mathcal{A}_{N}$
proceeds with the following modifications:

\begin{itemize}
\item $w=\varepsilon^{\beta}$ on $\partial\Omega$ replaced by $\frac{\partial
w}{\partial\nu}=0$ on $\partial\Omega$ in \eqref{pa};

\item $w=0$ on $\partial\Omega$ replaced by $\frac{\partial w}{\partial\nu}=0$
on $\partial\Omega$ in \eqref{lp};

\item no consideration of \eqref{pa2};

\item in Step 4, the test functions are now taken in $\mathcal{C}%
^{1}(\overline{\Omega})$.
\end{itemize}
\end{remark}

The following proposition shows that in general it is hard to give a lower
estimate for $q_{0}$ in Theorems \ref{t2} and \ref{t3}.

\begin{proposition}
\label{pex} \label{prop}Let $q\in\left(  0,1\right)  $. Then there exist
$\Omega$ and $a\in\mathcal{C}^{2}(\overline{\Omega})$ such that $q\not \in
\mathcal{A}_{D}$ and $q\not \in \mathcal{A}_{N}$.
\end{proposition}

\textbf{Proof}. Let $q\in\left(  0,1\right)  $. Define $\Omega:=\left(
0,\pi\right)  $,
\[
r:=\frac{2}{1-q}\in\left(  2,\infty\right)  ,\quad\text{and}\quad a\left(
x\right)  :=r^{1-\frac{2}{r}}\left(  1-r\cos^{2}x\right)  ,\quad\text{for
}x\in\overline{\Omega}.
\]
Clearly $a$ changes sign in $\Omega$. We now set%
\[
u\left(  x\right)  :=\frac{\sin^{r}x}{r}\in\mathcal{C}^{2}(\overline{\Omega
}).
\]
Note that $u>0$ in $\Omega$. We claim that
\[
\left\{
\begin{array}
[c]{lll}%
-u^{\prime\prime}=a(x)u^{q} & \mathrm{in} & \Omega,\\
u=u^{\prime}=0 & \mathrm{on} & \partial\Omega.
\end{array}
\right.
\]
Indeed, it is immediate to see that the boundary conditions are satisfied.
Also, taking into account that $rq=r-2$ (and so, $q=1-2/r$), a few
computations show that%
\begin{align*}
-u^{\prime\prime}  &  =-\left(  \left(  r-1\right)  \sin^{r-2}x\cos^{2}%
x-\sin^{r}x\right)  =\sin^{r-2}x\left(  1-r\cos^{2}x\right) \\
&  =\left(  1-r\cos^{2}x\right)  \sin^{rq}x=r^{1-\frac{2}{r}}\left(
1-r\cos^{2}x\right)  \left(  \frac{\sin^{r}x}{r}\right)  ^{q}=a\left(
x\right)  u^{q}%
\end{align*}
and therefore the claim follows.

To conclude the proof we note that, since $u>0$ in $\Omega$ and $u=u^{\prime
}=0$ on $\partial\Omega$, we have that $q\not \in \mathcal{A}_{D}$ and
$q\not \in \mathcal{A}_{N}$. $\blacksquare$

\begin{remark}
\label{core}Let $q$, $\Omega$, $a$ and $u$ be as in the above proposition.
Consider any bounded open interval $\Omega^{\prime}$ with $\Omega^{\prime
}\supset\overline{\Omega}$, and extend (to $\Omega^{\prime}$) $u$ by zero and
$a$ in any way. Then we clearly see that $u$ is a nontrivial nonnegative
solution having a \textit{dead core }in $\Omega^{\prime}$ (i.e. an open subset
with compact closure in $\Omega^{\prime}$ where $u$ vanishes) of both
$(P_{D})$ and $(P_{N})$, with $\Omega^{\prime}$ instead of $\Omega$.
\end{remark}

\section{Proofs of main results}

\begin{remark}
\label{rem:u0} The following fact shall be used several times in the sequel.
Let $\{u_{n}\}\subset H_{0}^{1}(\Omega)$ be a bounded sequence such that
%\marginpar{\textit{22oct2016}}%
\[
-\Delta u_{n}=h_{n}(x,u_{n})\quad\text{in }\Omega.
\]
Here $h_{n}:\Omega\times\left[  0,\infty\right)  \rightarrow\mathbb{R}$ are
Carath\'{e}odory functions satisfying
\[
|h_{n}(x,s)|\leq b(x)(1+s)\quad\mbox{ for }x\in\Omega,\mbox{ and }s\geq0,
\]
where $b\in L^{r}(\Omega)$, $r>N$. Then $\{u_{n}\}$ has a convergent
subsequence in $\mathcal{C}^{1}(\overline{\Omega})$. Indeed, by using the
above inequality on $h_{n}$, H\"{o}lder's inequality and the Sobolev embedding
theorem, we can derive that $\Vert u_{n}\Vert_{W^{2,\sigma_{k}}(\Omega)}$ is
bounded for each $\sigma_{k}=\frac{2N}{N-2k}$, $k=1,2,3,\dots$. Hence,
employing Sobolev's embedding theorem again, we obtain that $\Vert u_{n}%
\Vert_{\mathcal{C}^{1+\theta}(\overline{\Omega})}$ is bounded for some
$\theta\in(0,1)$. The desired conclusion follows by the Ascoli-Arzel\`{a}
theorem. We also note that a similar argument applies to the analogous Neumann problem.
\end{remark}

\textbf{Proof of Theorem \ref{t1}}. Assume by contradiction that $\left\{
a_{n}\right\}  $ is a sequence such that $a_{n}^{-}\rightarrow0$ in
$L^{r}(\Omega)$ and $u_{n}$ are nontrivial nonnegative solutions of $(P)$ with
$a=a_{n}$, satisfying that $u_{n}\not \in \mathcal{P}_{D}^{\circ}$. Let us
stress the fact that $a_{n}^{+}=a^{+}$ does not depend on $n$. We claim that
$\left\{  u_{n}\right\}  $ is bounded in $H_{0}^{1}(\Omega)$. Indeed, by our
assumptions on $f$, for any $\varepsilon>0$ there exists $C_{\varepsilon}>0$
such that
\begin{equation}
0\leq f(s)\leq C_{\varepsilon}+\varepsilon s,\quad\forall s\geq0. \label{cot}%
\end{equation}
Hence, for some $C,\tilde{C}_{\varepsilon}>0$, we have
\begin{align*}
\Vert u_{n}\Vert_{H_{0}^{1}(\Omega)}^{2}  &  =\int_{\Omega}|\nabla u_{n}%
|^{2}\leq\int_{\Omega}a^{+}(x)\left(  C_{\varepsilon}+\varepsilon
u_{n}\right)  u_{n}\\
&  \leq\tilde{C}_{\varepsilon}\Vert u_{n}\Vert_{H_{0}^{1}(\Omega
)}+C\varepsilon\Vert u_{n}\Vert_{H_{0}^{1}(\Omega)}^{2},
\end{align*}
where we have used Poincar\'{e}%
%TCIMACRO{\U{b4}}%
%BeginExpansion
\'{}%
%EndExpansion
s inequality. Taking $\varepsilon>0$ small enough, we deduce that $\left\{
u_{n}\right\}  $ is bounded in $H_{0}^{1}(\Omega)$. We can then assume that
$u_{n}\rightharpoonup u_{0}$ in $H_{0}^{1}(\Omega)$ and $u_{n}\rightarrow
u_{0}$ in $L^{p}(\Omega)$, with $p\in(1,2^{\ast})$, for some $u_{0}$. We claim
that $u_{0}\not \equiv 0$. Indeed, if $u_{0}\equiv0$ then, since
$u_{n}\rightarrow0$ in $L^{p}(\Omega)$, $a_{n}$ is bounded in $L^{r}(\Omega)$,
and
\begin{equation}
\int_{\Omega}\nabla u_{n}\nabla\phi=\int_{\Omega}a_{n}(x)f(u_{n})\phi
,\quad\forall\phi\in H_{0}^{1}(\Omega), \label{eun}%
\end{equation}
taking $\phi=u_{n}$ we see that $u_{n}\rightarrow0$ in $H_{0}^{1}(\Omega)$,
which contradicts Lemma \ref{l2}. Therefore $u_{0}\not \equiv 0$, as claimed.
In addition, since $u_{n}\rightharpoonup u_{0}$ in $H_{0}^{1}(\Omega)$,
recalling \eqref{cot} and choosing $\phi=u_{n}-u_{0}$ in \eqref{eun}, we
obtain $u_{n}\rightarrow u_{0}$ in $H_{0}^{1}(\Omega)$. Moreover, since
$a_{n}^{-}\rightarrow0$ in $L^{r}(\Omega)$, $u_{0}$ is a nontrivial
nonnegative solution of
\[
\left\{
\begin{array}
[c]{lll}%
-\Delta u_{0}=a^{+}(x)f(u_{0}) & \mathrm{in} & \Omega,\\
u_{0}=0 & \mathrm{on} & \partial\Omega.
\end{array}
\right.
\]
By the strong maximum principle and Hopf's Lemma, we have $u_{0}\in
\mathcal{P}_{D}^{\circ}$. Furthermore, standard elliptic regularity yields, up
to a subsequence, that $u_{n}\rightarrow u_{0}$ in $\mathcal{C}^{1}%
(\overline{\Omega})$ (see Remark \ref{rem:u0}, with $h_{n}(x,s)=a_{n}%
(x)f(s)$). Thus we must have $u_{n}\in\mathcal{P}_{D}^{\circ}$ for $n$ large
enough, which contradicts the assumption that $u_{n}\not \in \mathcal{P}%
_{D}^{\circ}$. $\blacksquare$\newline

\bigskip

\textbf{Proof of Theorem \ref{t2}}. First we note that, for every $c>0$, $u$
is a nonnegative solution of $(P_{D})$ if and only if $v:=c^{1/\left(
1-q\right)  }u$ is a nonnegative solution of $(P_{D})$ with $a$ replaced by
$ca$. Let $\Omega^{\prime}\neq\emptyset$ be a subdomain of $\Omega_{+}$. Since
$\lambda_{1}(ca,\Omega^{\prime})=c^{-1}\lambda_{1}(a,\Omega^{\prime
})\rightarrow0$ as $c\rightarrow\infty$, we can then assume without loss of
generality that $\lambda_{1}(a,\Omega^{\prime})<1$.

Assume by contradiction that $q_{n}\rightarrow1^{-}$ and $u_{n}$ are
nontrivial nonnegative solutions of $(P_{D})$ with $q=q_{n}$ and
$u_{n}\not \in \mathcal{P}_{D}^{\circ}$.

First we assume that $\{u_{n}\}$ is bounded in $H_{0}^{1}(\Omega)$. We can
assume that $u_{n}\rightharpoonup u_{0}$ in $H_{0}^{1}(\Omega)$,
$u_{n}\rightarrow u_{0}$ in $L^{p}(\Omega)$ with $p\in(1,2^{\ast})$, and
$u_{n}\rightarrow u_{0}$ \textit{a.e.} in $\Omega$, for some $u_{0}$. By Lemma
\ref{l3} (i), we have that $u_{0}\not \equiv 0$. Also, from
\[
\int_{\Omega}\nabla u_{n}\nabla(u_{n}-u_{0})=\int_{\Omega}a\left(  x\right)
u_{n}^{q_{n}}(u_{n}-u_{0})\rightarrow0,
\]
we infer that $u_{n}\rightarrow u_{0}$ in $H_{0}^{1}(\Omega)$. Moreover,
$u_{0}$ satisfies
\[
-\Delta u_{0}=a\left(  x\right)  u_{0},\quad u_{0}\geq0,\quad u_{0}\in
H_{0}^{1}(\Omega).
\]
We infer then, by Remark \ref{rem:u0} (with $h_{n}(x,s)=a(x)s^{q_{n}}$), that
$u_{n}\rightarrow u_{0}$ in $\mathcal{C}^{1}(\overline{\Omega})$, up to a
subsequence. Moreover, by the strong maximum principle and Hopf's Lemma we get
that $u_{0}\in\mathcal{P}_{D}^{\circ}$ and consequently $u_{n}\in
\mathcal{P}_{D}^{\circ}$ for $n$ large enough, which yields a contradiction.

We assume now that $\{u_{n}\}$ is unbounded in $H_{0}^{1}(\Omega)$. Then we
can assume that
\[
\left\Vert u_{n}\right\Vert :=\left\Vert u_{n}\right\Vert _{H_{0}^{1}(\Omega
)}\rightarrow\infty,\quad v_{n}:=\frac{u_{n}}{\left\Vert u_{n}\right\Vert
}\rightharpoonup v_{0}\text{ in }H_{0}^{1}(\Omega),
\]
and
\[
v_{n}\rightarrow v_{0}\text{ in }L^{p}(\Omega),\quad\mbox{with }p\in
(1,2^{\ast}),
\]
for some $v_{0}$. Note that $v_{n}$ satisfies
\begin{equation}
-\Delta v_{n}=a\left(  x\right)  \frac{v_{n}^{q_{n}}}{\left\Vert
u_{n}\right\Vert ^{1-q_{n}}},\quad v_{n}\geq0,\quad v_{n}\in H_{0}^{1}%
(\Omega). \label{ev}%
\end{equation}
Since $\left\Vert u_{n}\right\Vert \geq1$ for $n$ large enough, we have either
$\left\Vert u_{n}\right\Vert ^{1-q_{n}}\rightarrow\infty$ or $\left\Vert
u_{n}\right\Vert ^{1-q_{n}}$ is bounded. In the first case, from \eqref{ev} we
have
\[
1=\left\Vert v_{n}\right\Vert ^{2}=\frac{\int_{\Omega}a\left(  x\right)
v_{n}^{q_{n}+1}}{\left\Vert u_{n}\right\Vert ^{1-q_{n}}}\rightarrow0,
\]
which is a contradiction. Now, if $\left\Vert u_{n}\right\Vert ^{1-q_{n}}$ is
bounded then we can assume that $\left\Vert u_{n}\right\Vert ^{1-q_{n}%
}\rightarrow d\geq1$. From \eqref{ev}, we obtain
\[
\int_{\Omega}\nabla v_{0}\nabla\phi=\frac{1}{d}\int_{\Omega}a\left(  x\right)
v_{0}\phi,\quad\forall\phi\in H_{0}^{1}(\Omega),
\]
i.e.
\[
-\Delta v_{0}=\frac{1}{d}a\left(  x\right)  v_{0}\quad\text{in }\Omega,\quad
v_{0}\in H_{0}^{1}(\Omega).
\]
In addition, $v_{n}\rightarrow v_{0}$ in $H_{0}^{1}(\Omega)$, so that
$v_{0}\not \equiv 0$ and $v_{0}\geq0$. Once again, by the strong maximum
principle and Hopf's Lemma, we deduce that $v_{0}\in\mathcal{P}_{D}^{\circ}$.
Furthermore, recalling Remark \ref{rem:u0} (with $h_{n}(x,s)=a(x)\Vert
u_{n}\Vert^{q_{n}-1}s^{q_{n}}$) we have that $v_{n}\rightarrow v_{0}$ in
$\mathcal{C}^{1}(\overline{\Omega})$, up to a subsequence. Consequently
$v_{n}\in\mathcal{P}_{D}^{\circ}$ for $n$ large enough. Hence $u_{n}%
\in\mathcal{P}_{D}^{\circ}$, and we get another contradiction, which concludes
the proof. $\blacksquare$\newline

\bigskip

\textbf{Proof of Theorem \ref{t3}}. We proceed as in the proof of Theorem
\ref{t2}: we assume by contradiction that $q_{n}\rightarrow1^{-}$ and $u_{n}$
are nontrivial nonnegative solutions of $(P_{N})$ with $q=q_{n}$, and that
$u_{n}$ do not belong to $\mathcal{P}_{N}^{\circ}$.

First we suppose that $\{u_{n}\}$ is unbounded in $H^{1}(\Omega)$. Then we can
assume that
\[
\left\Vert u_{n}\right\Vert :=\left\Vert u_{n}\right\Vert _{H^{1}(\Omega
)}\rightarrow\infty,\quad v_{n}:=\frac{u_{n}}{\left\Vert u_{n}\right\Vert
}\rightharpoonup v_{0}\text{ in }H^{1}(\Omega),
\]
and
\[
v_{n}\rightarrow v_{0}\text{ in }L^{p}(\Omega),\quad\mbox{with }p\in
(1,2^{\ast}),
\]
for some $v_{0}$. Note that $v_{n}$ satisfies
\begin{equation}
-\Delta v_{n}=a\left(  x\right)  \frac{v_{n}^{q_{n}}}{\left\Vert
u_{n}\right\Vert ^{1-q_{n}}},\quad v_{n}\geq0,\quad v_{n}\in H^{1}%
(\Omega),\label{evn}%
\end{equation}
and so we have
\[
\int_{\Omega}\nabla v_{n}\nabla\phi=\int_{\Omega}a(x)\frac{v_{n}^{q_{n}}%
}{\left\Vert u_{n}\right\Vert ^{1-q_{n}}}\phi,\quad\forall\phi\in H^{1}%
(\Omega).
\]
Now, if $\Vert u_{n}\Vert^{1-q_{n}}\rightarrow\infty$, taking $\phi=v_{n}$ we
obtain that $\int_{\Omega}|\nabla v_{n}|^{2}\rightarrow0$, which implies that
$v_{n}\rightarrow v_{0}$ in $H^{1}(\Omega)$ and $v_{0}$ is a nonnegative
constant. Since $\left\Vert v_{n}\right\Vert =1$, we infer that $v_{0}$ is a
positive constant. By Remark \ref{rem:u0}, we have $v_{n}\rightarrow v_{0}$ in
$\mathcal{C}^{1}(\overline{\Omega})$. Consequently, $v_{n}\in\mathcal{P}%
_{N}^{\circ}$ for $n$ large enough, which yields a contradiction.

On the other hand, if $\left\Vert u_{n}\right\Vert ^{1-q_{n}}$ is bounded,
reasoning as in the proof of Theorem \ref{t2} we derive now that $v_{0}$
satisfies
\[
-\Delta v_{0}=\frac{1}{d}a\left(  x\right)  v_{0}\quad\text{in }\Omega
,\quad\text{with}\quad\frac{\partial v_{0}}{\partial\nu}=0\text{ on }%
\partial\Omega.
\]
In addition, $v_{n}\rightarrow v_{0}$ in $H^{1}(\Omega)$, so that
$v_{0}\not \equiv 0$ and $v_{0}\geq0$. By the strong maximum principle, we
deduce that $v_{0}\in\mathcal{P}_{N}^{\circ}$. Once again, by Remark
\ref{rem:u0}, we have $v_{n}\rightarrow v_{0}$ in $\mathcal{C}^{1}%
(\overline{\Omega})$, so that again we reach a contradiction.

Finally, if $\{u_{n}\}$ is bounded in $H^{1}(\Omega)$ then we can argue again
as in the proof of Theorem \ref{t2}. Indeed, by Lemma \ref{l3} (ii), we have
that $v_{0}\not \equiv 0$. The rest of the proof is similar, so we omit it.
$\blacksquare$\newline

\textbf{Proof of Corollary \ref{c1}}. We first claim that $(P)$ admits a
nontrivial nonnegative solution $u$. Indeed, let $\underline{u}_{\gamma
}:=\varphi_{\gamma,0}$, where $\varphi_{\gamma,0}$ is given by (\ref{zzz})
(with $\delta=\gamma$ and $\varepsilon=0$). Using the condition that
$\displaystyle\lim_{s\rightarrow0^{+}}\frac{f(s)}{s}=\infty$ and arguing as in
the first part of the proof of Lemma \ref{l1}, it is easy to see that for all
$\gamma>0$ sufficiently small, $\underline{u}_{\gamma}$ is a (nonnegative)
subsolution of $(P)$. On the other side, let $\psi>0$ be the unique solution
of the problem%
\[
\left\{
\begin{array}
[c]{lll}%
-\Delta\psi=a^{+}(x) & \mathrm{in} & \Omega,\\
\psi=0 & \mathrm{on} & \partial\Omega.
\end{array}
\right.
\]
By utilizing $\displaystyle\lim_{s\rightarrow\infty}\frac{f(s)}{s}=0$, we note
that for every $k>0$ large enough, $\overline{u}_{k}:=k\psi$ is a
supersolution of $(P)$. Indeed, there exists $C>0$ such that
\[
f(s)\leq\frac{1}{2\Vert\psi\Vert_{L^{\infty}(\Omega)}}s+C,\quad s\geq0.
\]
It follows that
\[
-\Delta(k\psi)-a(x)f(k\psi)\geq a^{+}(x)\left(  \frac{k}{2}-C\right)
\geq0\quad\mbox{in }\Omega,\quad\mbox{ if }k\geq2C.
\]
Moreover, since $\underline{u}_{\gamma}=0$ in a neighborhood of $\partial
\Omega$, making $\delta$ smaller and $k$ larger if necessary, we have that
$\underline{u}_{\gamma}\leq\overline{u}_{k}$, and it follows that $(P)$ has a
nontrivial nonnegative solution $u$.

Since $\left\Vert a^{-}\right\Vert _{L^{r}(\Omega)}<\delta$, from Theorem
\ref{t1}, we know that any nontrivial nonnegative solution of $(P)$ belongs to
$\mathcal{P}_{D}^{\circ}$. In particular, $u\in\mathcal{P}_{D}^{\circ}$.

Finally, by the assumptions on $f$ and Theorem 2.1 in \cite{DS}, we also know
that there is at most one positive solution of $(P)$. Therefore there are no
other nontrivial nonnegative solutions of $(P)$. $\blacksquare$\newline

\textbf{Proof of Corollary \ref{c2}}. The proof is similar to the previous
one. It suffices to note that $f(s)=s^{q}$ satisfies $(H_{1})$, so that the
previous proof and Theorem \ref{t2} yield the existence assertion. In
addition, $f(s)=s^{q}$ satisfies the conditions of Corollary \ref{c1}, so that
the nonexistence assertion is proved in the same way. $\blacksquare$\newline

\textbf{Proof of Corollary \ref{c3}}. We set now
\[
I(u):=\int_{\Omega}\left(  \frac{1}{2}|\nabla u|^{2}-\frac{1}{q+1}%
a(x)|u|^{q+1}\right)  ,\quad\text{for }u\in H^{1}(\Omega).
\]
We claim that $I$ is coercive on $H^{1}(\Omega)$. Indeed, assume that
\[
u_{n}\in H^{1}(\Omega),\quad\left\Vert u_{n}\right\Vert :=\left\Vert
u_{n}\right\Vert _{H^{1}(\Omega)}\rightarrow\infty,\quad\text{and}\quad
I(u_{n})\text{ is bounded from above}.
\]
We may assume that
\[
v_{n}:=\frac{u_{n}}{\left\Vert u_{n}\right\Vert }\rightharpoonup v_{0}\text{
in }H^{1}(\Omega),\quad\text{and}\quad v_{n}\rightarrow v_{0}\text{ in }%
L^{p}(\Omega)\text{ for }p\in(1,2^{\ast}),
\]
for some $v_{0}$. Then
\[
\frac{1}{2}\int_{\Omega}|\nabla v_{n}|^{2}-\frac{1}{(q+1)\left\Vert
u_{n}\right\Vert ^{1-q}}\int_{\Omega}a(x)|v_{n}|^{q+1}=\frac{I(u_{n}%
)}{\left\Vert u_{n}\right\Vert ^{2}},
\]
so that $\int_{\Omega}|\nabla v_{n}|^{2}\rightarrow0$. It follows that
$v_{n}\rightarrow v_{0}$ in $H^{1}(\Omega)$ and $v_{0}$ is a nonzero constant.
Moreover, from
\[
-\frac{1}{q+1}\int_{\Omega}a(x)|v_{n}|^{q+1}<\frac{I(u_{n})}{\left\Vert
u_{n}\right\Vert ^{q+1}},
\]
we have that
\[
\int_{\Omega}a(x)|v_{0}|^{q+1}=\lim\int_{\Omega}a(x)|v_{n}|^{q+1}\geq0,
\]
and consequently $\int_{\Omega}a\geq0$, which contradicts our assumption.
Therefore $I$ is coercive so that it has a global maximum. Taking $u_{0}$ such
that $\int_{\Omega}a(x)|u_{0}|^{q+1}>0$, we see that $I(tu_{0})<0$ if $t>0$ is
sufficiently small. This shows that $I$ has a nontrivial global minimizer.
Finally, since $I$ is even, it has a nonnegative global minimizer, which is a
nontrivial nonnegative solution of $(P^{\prime})$. By Theorem \ref{t3}, this
solution (and any other nontrivial nonnegative solution) belongs to
$\mathcal{P}_{N}^{\circ}$ for $q_{0}<q<1$.

Lastly, reasoning exactly as in Lemma 3.1 in \cite{BPT}, we infer that there
are no other nontrivial nonnegative solutions of $(P_{N})$. $\blacksquare
$\newline

\textbf{Proof of Theorem \ref{ta}.} Note that Theorem \ref{t2} says that
$\mathcal{A}_{D}$ is nonempty. Now, first we show, via the continuity argument
used in the proofs of Theorems \ref{t2} and \ref{t3}, that $\mathcal{A}_{D}$
is open. Indeed, assume to the contrary that there exist $q\in\mathcal{A}_{D}$
and $q_{n}\not \in \mathcal{A}_{D}$ such that $q_{n}\rightarrow q$. We take
nontrivial nonnegative solutions $u_{n}\not \in \mathcal{P}_{D}^{\circ}$ of
$(P_{D})$ with $q=q_{n}$. Using Lemma \ref{l3} and arguing as in the proof of
Theorem \ref{t2}, we may deduce that $\{u_{n}\}$ is bounded in $H_{0}%
^{1}(\Omega)$ and consequently, up to a subsequence, $u_{n}\rightarrow u_{0}$
in $\mathcal{C}^{1}(\overline{\Omega})$, where $u_{0}$ is a nontrivial
nonnegative solution of $(P_{D})$. Since $q\in\mathcal{A}_{D}$, we have
$u_{0}\in\mathcal{P}_{D}^{\circ}$, and so $u_{n}\in\mathcal{P}_{D}^{\circ}$
for $n$ large enough, which is a contradiction. Thus $\mathcal{A}_{D}$ is open.

We prove next that $\mathcal{A}_{D}$ is connected. To this end, we show that
if $q_{0}\in\mathcal{A}_{D}$ then $(q_{0},1)\subset\mathcal{A}_{D}$. Let
$q_{0}\in\mathcal{A}_{D}$. By Lemma \ref{l6}, we know that
\[
\left(  q_{0},\frac{1}{2-q_{0}}-\sigma_{0}\right]  \subset\mathcal{A}_{D}%
\]
with $\sigma_{0}=\frac{1}{10}\frac{(q_{0}-1)^{2}}{2-q_{0}}$, where
$\frac{(q_{0}-1)^{2}}{2-q_{0}}$ is the length of the interval $(q_{0},\frac
{1}{2-q_{0}})$. By iteration, since $q_{1}=\frac{1}{2-q_{0}}-\sigma_{0}\in A$,
we have, again by Lemma \ref{l6}, that
\[
\left(  q_{1},\frac{1}{2-q_{1}}-\sigma_{1}\right]  \subset\mathcal{A}_{D},
\]
where $\sigma_{1}=\frac{1}{10}\frac{(q_{1}-1)^{2}}{2-q_{1}}$. More generally,
we have
\[
\left(  q_{n-1},\frac{1}{2-q_{n-1}}-\sigma_{n-1}\right]  \subset
\mathcal{A}_{D},
\]
where $\sigma_{n}:=\frac{1}{10}\frac{(q_{n}-1)^{2}}{2-q_{n}}$, and
$q_{n}:=\frac{1}{2-q_{n-1}}-\sigma_{n-1}$. Then, we obtain by induction that
$\{q_{n}\}$ is nondecreasing and $q_{n}\leq1$, so that $q_{n}\rightarrow
q_{\ast}$ for some $q_{\ast}\leq1$. Passing to the limit as $n\rightarrow
\infty$, we have
\[
q_{\ast}=\frac{1}{2-q_{\ast}}-\frac{1}{10}\frac{(q_{\ast}-1)^{2}}{2-q_{\ast}},
\]
so that $\frac{9}{10}\frac{(q_{\ast}-1)^{2}}{2-q_{\ast}}=0$, and thus,
$q_{\ast}=1$. Hence we have proved that $(q_{0},1)\subset\mathcal{A}_{D}$.

Finally, the proof that $\mathcal{A}_{N}$ is open and connected can be carried
out in the same manner as for $\mathcal{A}_{D}$. In addition, we know by
Corollary \ref{c3} that $\int_{\Omega}a<0$ is sufficient for the existence of
some $q\in\mathcal{A}_{N}$. The proof is now complete. $\blacksquare$\newline

\section{Positivity results for concave-convex type problems}

%Throughout this section we fix $f\left(  s\right)  =s^{q}$.
As an application of some of our previous results, we consider now the
problem
\[
\left\{
\begin{array}
[c]{lll}%
-\Delta u=\lambda a(x)u^{q}+g(u) & \mathrm{in} & \Omega,\\
u=0 & \mathrm{on} & \partial\Omega,
\end{array}
\right.  \leqno{(P_\lambda)}
\]
where now $a\in L^{\infty}(\Omega)$, $0<q<1$, $\lambda>0$, and $N\geq3$. In
addition, we assume that $g:[0,\infty)\rightarrow\lbrack0,\infty)$ is
continuous and superlinear in the following sense:
\[
\lim_{s\rightarrow0^{+}}\frac{g(s)}{s}=0\quad\text{and}\quad\lim
_{s\rightarrow\infty}\frac{g(s)}{s^{p}}=1,\quad\text{for some }1<p<\frac
{N+2}{N-2}.\leqno{(H_4)}
\]
We shall also assume that $g(s)>0$ for $s>0$.

The problem above has been investigated in \cite{ABC} for $a\equiv1$, and
$g(s)=s^{p}$. The authors proved that $(P_{\lambda})$ has two positive
solutions for $\lambda>0$ sufficiently small. This result was extended to a
more general nonlinearity, with $a\geq0$, in \cite{DGU2}. In addition, in
\cite{DGU1}, the authors allowed $a$ to change sign and proved the existence
of two nontrivial nonnegative solutions of $(P_{\lambda})$.

The growth condition at infinity in $(H_{4})$ ensures, in particular, an
\textit{a priori} bound for nonnegative solutions of $(P_{\lambda})$ (see
\cite{ALG98}, and also \cite{LGMMT13}), which will be used to prove the
following positivity result:

\begin{theorem}
\label{t4} Assume $(H_{2})$ and $(H_{4})$. In addition, assume that every
nontrivial nonnegative solution of $(P_{D})$ belongs to $\mathcal{P}%
_{D}^{\circ}$. Then there exists $\lambda_{0}>0$ such that every nontrivial
nonnegative solution of $(P_{\lambda})$ belongs to $\mathcal{P}_{D}^{\circ}$
for $0<\lambda<\lambda_{0}$.
\end{theorem}

\textbf{Proof}. Assume by contradiction that $\lambda_{n}\rightarrow0^{+}$ and
$u_{n}$ are nontrivial nonnegative solutions of $(P_{\lambda})$ with
$\lambda=\lambda_{n}$ and $u_{n}\not \in \mathcal{P}_{D}^{\circ}$. By the
\textit{a priori} bounds in \cite{ALG98}, there exists $K>0$ such that $\Vert
u_{n}\Vert_{\infty}\leq K$ for every $n$. It follows that $\{u_{n}\}$ is
bounded in $H_{0}^{1}(\Omega)$. Thus, we can assume that $u_{n}\rightharpoonup
u_{0}$ in $H_{0}^{1}(\Omega)$ and $u_{n}\rightarrow u_{0}$ in $L^{s}(\Omega)$,
$1<s<2^{\ast}$, for some $u_{0}$. Taking as test function $u_{n}-u_{0}$ in
$(P_{\lambda_{n}})$ we see that $u_{n}\rightarrow u_{0}$ in $H_{0}^{1}%
(\Omega)$ and $u_{0}$ is a solution of $(P_{\lambda})$ with $\lambda=0$.
Moreover, by elliptic regularity, we have, up to a subsequence, that
$u_{n}\rightarrow u_{0}$ in $\mathcal{C}^{1}(\overline{\Omega})$, see Remark
\ref{rem:u00} below. If $u_{0}\not \equiv 0$ then, by the strong maximum
principle, we have that $u_{0}\in\mathcal{P}_{D}^{\circ}$, and consequently
$u_{n}\in\mathcal{P}_{D}^{\circ}$ for $n$ large enough, which provides a
contradiction. Now, if $u_{0}\equiv0$, then we consider $v_{n}:=\lambda
_{n}^{\frac{1}{q-1}}u_{n}$. We see that $v_{n}$ are nontrivial nonnegative
solutions of
\[
\left\{
\begin{array}
[c]{lll}%
-\Delta v=a(x)v^{q}+\lambda^{\frac{1}{q-1}}g(\lambda^{\frac{1}{1-q}}v) &
\mathrm{in} & \Omega,\\
v=0 & \mathrm{on} & \partial\Omega,
\end{array}
\right.
\]
with $\lambda=\lambda_{n}$. Hence, $v_{n}$ are nontrivial nonnegative
supersolutions of
\[
\left\{
\begin{array}
[c]{lll}%
-\Delta v=a(x)v^{q} & \mathrm{in} & \Omega,\\
v=0 & \mathrm{on} & \partial\Omega,
\end{array}
\right.
\]
with $\lambda=\lambda_{n}$.

We claim that $v_{n}\not \equiv 0$ in $\Omega_{+}$. Indeed, note that if
$v_{n}\equiv0$ in $\Omega_{+}$ then $u_{n}\equiv0$ in $\Omega_{+}$, so that
\[
\int_{\Omega}|\nabla u_{n}|^{2}\leq\int_{\Omega}g(u_{n})u_{n}\leq
\varepsilon\Vert u_{n}\Vert_{H_{0}^{1}(\Omega)}^{2}+C_{\varepsilon}\Vert
u_{n}\Vert_{H_{0}^{1}(\Omega)}^{p+1}.
\]
Taking $\varepsilon>0$ sufficiently small we see that $\Vert u_{n}\Vert
_{H_{0}^{1}(\Omega)}\geq C>0$, which contradicts $u_{n}\rightarrow0$ in
$H_{0}^{1}(\Omega)$. Therefore the claim is proved. Since $\Omega_{+}$ has
finitely many connected components, we can assume that $v_{n}\not \equiv 0$ in
some fixed subdomain $\Omega^{\prime}\subset\Omega_{+}$. Let $\phi$ be as in
the proof of Lemma \ref{l1}. Arguing as in this proof, we have that
$\varepsilon\phi$ is a nonnegative subsolution of
\[
-\Delta u=a(x)u^{q}\text{\quad in }B,
\]
where $B$ is an open ball such that $\overline{B}\subset\Omega^{\prime}$. We
extend $\phi$ by zero to $\overline{\Omega}\setminus\overline{B}$. For
$\varepsilon>0$ small enough, we have that $\varepsilon\phi\leq v_{n}$ for
every $n$. Thus, we find a nonnegative solution $w_{n}$ of $(P)$ such that
$\varepsilon\phi\leq w_{n}\leq v_{n}$. But, by our assumption, we have that
$w_{n}\in\mathcal{P}_{D}^{\circ}$, which contradicts the assumption that
$v_{n}\not \in \mathcal{P}_{D}^{\circ}$. The proof is now complete.
$\blacksquare$\newline

\begin{remark}
\label{rem:u00} In the same way as Remark \ref{rem:u0}, we give some further
details on the regularity argument used in the previous proof. We set now
\[
h_{\lambda}(x,s):=\lambda a(x)f(s)+g(s)
\]
and use the conditions
\[
a\in L^{\infty}(\Omega),\text{ \ }\lim_{s\rightarrow\infty}\frac{f(s)}%
{s}=0,\text{ \ }\mbox{ and }\text{ \ }\lim_{s\rightarrow\infty}\frac
{g(s)}{s^{p}}=1\text{ \ }\mbox{ for some }\text{ \ }1<p<\frac{N+2}{N-2},
\]
to infer that, given $\overline{\lambda}>0$, there exists $C>0$ such that
\[
|h_{\lambda}(x,s)|\leq C(1+s^{p})\ \ \mbox{ for }x\in\Omega,\ s\geq
0,\mbox{ and }|\lambda|\leq\overline{\lambda}.
\]
In the same manner as in Remark \ref{rem:u0}, we can deduce that $\Vert
u_{n}\Vert_{W^{2,\sigma_{k}}(\Omega)}$ is bounded for each $\sigma_{k}%
=\frac{2N}{p_{k}}$, $k=1,2,3,\dots$, where
\[
p_{k}=\frac{4p}{p-1}-p^{k}\left(  \frac{4}{p-1}-(N-2)\right)  .
\]
Since $p>1$ and $\frac{4}{p-1}-(N-2)>0$, we can choose $\sigma_{k}>N$ such
that $\Vert u_{n}\Vert_{W^{2,\sigma_{k}}(\Omega)}$ is bounded. Then, the
argument proceeds in the same way as in Remark \ref{rem:u0}.
\end{remark}

As a consequence of Theorem \ref{t4}, we obtain two positive solutions of
$(P_{\lambda})$ for $\lambda>0$ small, if either $a^{-}$ is small or $q$ is
close to $1$:

\begin{corollary}
\label{c4} Assume $(H_{2})$ and $(H_{4})$. Then there exist $\delta>0$ and
$q_{0}\in(0,1)$ such that, if either $\Vert a^{-}\Vert_{L^{r}\left(
\Omega\right)  }<\delta$ or $q_{0}<q<1$, then there exists $\lambda_{0}>0$
with the following properties:

\begin{enumerate}
\item[(i)] any nontrivial nonnegative solution of $(P_{\lambda})$ belongs to
$\mathcal{P}_{D}^{\circ}$ for $0<\lambda<\lambda_{0}$.

\item[(ii)] $(P_{\lambda})$ has two solutions in $\mathcal{P}^{\circ}_{D}$ for
$0<\lambda<\lambda_{0}$.\newline
\end{enumerate}
\end{corollary}

\textbf{Proof.}

\begin{enumerate}
\item[(i)] We apply Theorems \ref{t1} and \ref{t2} to $(P_{D})$ and obtain,
respectively, $\delta>0$ and $q_{0}\in(0,1)$ such that every nontrivial
nonnegative solution of $(P_{D})$ belongs to $\mathcal{P}_{D}^{\circ}$ if
either $\Vert a^{-}\Vert_{L^{r}\left(  \Omega\right)  }<\delta$ or $q_{0}%
<q<1$. Theorem \ref{t4} yields the conclusion.

\item[(ii)] We use Theorem 2.1 from \cite{DGU1}. One can easily show that
assumptions $(H_{0})$-$(H_{5})$ from \cite{DGU1} are satisfied under our
conditions. Thus there exists $\lambda>0$ such that for $0<\lambda<\lambda
_{0}$ there exist two nontrivial nonnegative solutions of $(P_{\lambda})$.
Decreasing $\lambda_{0}$ if necessary, by the previous item, we infer that
these solutions belong to $\mathcal{P}_{D}^{\circ}$ if either $\Vert
a^{-}\Vert_{L^{r}\left(  \Omega\right)  }<\delta$ or $q_{0}<q<1.$
$\blacksquare$\newline
\end{enumerate}

\begin{remark}
Let us set
\[
\mathcal{B}_{D}:=\{\lambda
>0:\mbox{any nontrivial non-negative solution of $(P_\lambda)$ belongs to $\mathcal{P}^\circ_{D}$}\},
\]
and assume $(H_{2})$, $(H_{4})$, and either $\Vert a^{-}\Vert_{L^{r}\left(
\Omega\right)  }<\delta$ or $q_{0}<q<1$, where $\delta$ and $q_{0}$ are
provided by Corollary \ref{c4}. Then $(0,\lambda_{0})\subset\mathcal{B}_{D}$.
Moreover, arguing as in the proof of Theorem \ref{t4} we can show that
$\mathcal{B}_{D}$ is an open set. Indeed, one may easily see that the proof of
Theorem \ref{t4} carries on taking now a sequence $\lambda_{n}\rightarrow
\lambda_{0}\in\mathcal{B}_{D}$. 
\end{remark}

We establish now a result analogue to Theorem \ref{t4} for the problem
\[
\left\{
\begin{array}
[c]{lll}%
-\Delta u=\lambda a(x)u^{q}+g(u) & \mathrm{in} & \Omega,\\
\frac{\partial u}{\partial\nu}=0 & \mathrm{on} & \partial\Omega,
\end{array}
\right.  \leqno{(Q_\lambda)}
\]
Instead of $(H_{4})$, we shall assume now
\[
\lim_{s\rightarrow0^{+}}\frac{g(s)}{s}=0\quad\text{and}\quad\lim
_{s\rightarrow\infty}\frac{g(s)}{s^{p}}=1,\text{ for some }1<p<\frac{N+1}%
{N-1}.\leqno{(H_{4}')}
\]

The above problem, with $g(s)=s^{p}$,  has been recently investigated in
\cite{RQU6}. The authors established existence and multiplicity results for
nontrivial nonnegative solutions of $(Q_{\lambda})$, for $\lambda>0$
sufficiently small. Furthermore, the asymptotic behavior of these solutions as
$\lambda\rightarrow0^{+}$ provides the positivity of some of these solutions,
in certain cases. We shall now prove a general positivity result for
$(Q_{\lambda})$:
%\marginpar{\textit{22oct2016}}

\begin{theorem}
\label{t5} Assume $(H_{2})$ and $(H_{4}^{\prime})$. In addition, assume that
every nontrivial nonnegative solution of $(P_{N})$ belongs to $\mathcal{P}%
_{N}^{\circ}$. Then there exists $\lambda_{0}>0$ such that every nontrivial
nonnegative solution of $(Q_{\lambda})$ belongs to $\mathcal{P}_{N}^{\circ}$
for $0<\lambda<\lambda_{0}$.
\end{theorem}

\textbf{Proof}. The proof is similar to the one of Theorem \ref{t4}. Again by
the \textit{a priori }bounds of \cite{ALG98} we get that $\{u_{n}\}$ is
bounded in $H^{1}(\Omega)$. To show that if $u_{n}\rightarrow0$ in
$H^{1}(\Omega)$ then $u_{n}\not \equiv 0$ in $\Omega_{+}$, we proceed in the
following way: assume the contrary and set $w_{n}:=\frac{u_{n}}{\Vert
u_{n}\Vert}$, where $\Vert u_{n}\Vert:=\Vert u_{n}\Vert_{H^{1}(\Omega)}$. So
we can assume that $w_{n}\rightharpoonup w_{0}$ in $H^{1}(\Omega)$,
$w_{n}\rightarrow w_{0}$ in $L^{s}(\Omega)$, with $1<s<2^{\ast}$, and
$w_{n}\rightarrow w_{0}$ \textit{a.e.} in $\Omega$, for some $w_{0}$. Thus,
from
\[
\int_{\Omega}|\nabla u_{n}|^{2}\leq\int_{\Omega}g(u_{n})u_{n},
\]
we obtain that
\[
\int_{\Omega}|\nabla w_{n}|^{2}\leq\int_{\Omega}\frac{g(\Vert u_{n}\Vert
w_{n})}{\Vert u_{n}\Vert}w_{n}=\int_{\text{supp }w_{n}}\frac{g(\Vert
u_{n}\Vert w_{n})}{\Vert u_{n}\Vert w_{n}}w_{n}^{2}%
\]
Since $\Vert u_{n}\Vert\rightarrow0$ and $\frac{g(s)}{s}\rightarrow0$ as
$s\rightarrow0^{+}$, we easily see that
\[
\int_{\Omega}|\nabla w_{n}|^{2}\rightarrow0,
\]
so that $w_{n}\rightarrow w_{0}$ in $H^{1}(\Omega)$ and $w_{0}$ is a positive
constant. This contradicts the assumption that $w_{n}\equiv0$ in $\Omega_{+}$.
The rest of the proof carries on in a similar way. $\blacksquare$

\begin{corollary}
\label{c5} Assume $(H_{2})$ and $(H_{4}^{\prime})$. Let $q_{0}\in(0,1)$ be
given by Theorem \ref{t3}. If $q_{0}<q<1$ then there exists $\lambda_{0}>0$
such that any nontrivial nonnegative solution of $(Q_{\lambda})$ belongs to
$\mathcal{P}_{N}^{\circ}$ for $0<\lambda<\lambda_{0}$.
\end{corollary}

\begin{corollary}
\label{c6} Let $g(s)=s^{p}$, with $1<p<\frac{N+1}{N-1}$, and assume
$q_{0}<q<1$, where $q_{0}\in(0,1)$ is given by Theorem \ref{t3}. If
$\int_{\Omega}a<0$ then there exists $\lambda_{0}>0$ such that $(Q_{\lambda})$
has two solutions in $\mathcal{P}^{\circ}_{N}$ for $0<\lambda<\lambda_{0}$.
\end{corollary}

\textbf{Proof}. We apply Corollary 1.3 (2) from \cite{RQU6} to obtain
$\lambda_{0}>0$ such that $(Q_{\lambda})$ has two solutions $u_{1,\lambda}$,
$u_{2,\lambda}$ such that $u_{2,\lambda}>u_{1,\lambda}\geq0$ in $\overline
{\Omega}$ for $0<\lambda<\lambda_{0}$. By Corollary \ref{c5}, decreasing
$\lambda_{0}$ if necessary, we have that $u_{1,\lambda}$ and $u_{2,\lambda}$
belong to $\mathcal{P}_{N}^{\circ}$. $\blacksquare$

\begin{remark}
Let us set
\[
\mathcal{B}_{N}:=\{\lambda
>0:\mbox{any nontrivial non-negative solution of $(Q_\lambda)$ belongs to $\mathcal{P}^\circ_{N}$}\}.
\]
We assume that $(H_{2})$, $(H_{4}^{\prime})$ hold, and $q_{0}<q<1$, where
$q_{0}$ is provided by Corollary \ref{c5}. Then, by arguing in the same way as
for $(P_{\lambda})$, we observe that $(0, \lambda_{0}) \subset\mathcal{B}_{N}%
$, and in addition, $\mathcal{B}_{N}$ is open.
%$(0,\lambda_{0})\subset\mathcal{B}_{N}$. Moreover, arguing as in the proof of Theorem \ref{t5} we can show that $\mathcal{B}_{N}$ is an open set.
%Indeed, one may easily see that the proof of Theorem \ref{t4} carries on taking now a sequence $\lambda_{n}\rightarrow \lambda_{0}\in\mathcal{B}_{D}$. $\blacksquare$

\end{remark}

%%%%
\subsection*{Acknowledgements}

U. Kaufmann was partially supported by Secyt-UNC. H. Ramos Quoirin was supported by Fondecyt 1161635. K. Umezu was supported by JSPS KAKENHI Grant Number 15K04945.


\begin{thebibliography}{99}                                                                                               %


\bibitem {alama}S. Alama, \textit{Semilinear elliptic equations with sublinear
indefinite nonlinearities}, Adv. Differential Equations \textbf{4} (1999), 813--842.

\bibitem {ALG98}H. Amann, J. L\'{o}pez-G\'{o}mez, \textit{A priori bounds and
multiple solutions for superlinear indefinite elliptic problems},
J.\ Differential Equations \textbf{146} (1998), 336--374.

\bibitem {ABC}A. Ambrosetti, H. Brezis, G. Cerami, \textit{Combined effects of
concave and convex nonlinearities in some elliptic problems}, J. Funct. Anal.
\textbf{122} (1994), 519-543.

\bibitem {bandle}C. Bandle, M. Pozio, A. Tesei, \textit{The asymptotic
behavior of the solutions of degenerate parabolic equations,} Trans. Amer.
Math. Soc. \textbf{303} (1987), 487-501.\textit{\ }

\bibitem {BPT}C.\ Bandle, A.\ M.\ Pozio, A.\ Tesei, \textit{Existence and
uniqueness of solutions of nonlinear Neumann problems},
Math.\ Z.\ \textbf{199} (1988), 257--278.

%\bibitem {d}D. G. de Figueiredo, Positive solutions of semilinear elliptic
%equations, Lect. Notes Math. Springer \textbf{957}, 34--87 (1982).


\bibitem {cuesta}M. Cuesta, P. Tak\'{a}\v{c}, \textit{A strong comparison
principle for positive solutions of degenerate elliptic equations},
Differential Integral Equations \textbf{13} (2000), 721--746.

\bibitem {DFG}D. G. de Figueiredo, J.-P. Gossez, \textit{On the first curve of
the Fucik spectrum of an elliptic operator}, Differential Integral Equations
\textbf{7} (1994), 1285-1302.

\bibitem {DGU1}D. G. de Figueiredo, J-P. Gossez, P. Ubilla, \textit{Local
superlinearity and sublinearity for indefinite semilinear elliptic problems},
J. Funct. Anal. \textbf{199} (2003), 452--467.

\bibitem {DGU2}D. G. De Figueiredo, J-P. Gossez, P. Ubilla,
\textit{Multiplicity results for a family of semilinear elliptic problems
under local superlinearity and sublinearity}, J. Eur. Math. Soc. \textbf{8}
(2006), 269--286.

\bibitem {DS}M. Delgado, A. Su\'{a}rez, \textit{On the uniqueness of positive
solution of an elliptic equation}, Appl. Math. Lett. \textbf{18} (2005), 1089-1093.

\bibitem {du}Y. Du, \textit{Order structure and topological methods in
nonlinear partial differential equations. Vol. 1. Maximum principles and
applications}, World Scientific Publishing Co. Pte. Ltd., Hackensack, NJ, 2006.

\bibitem {nodea}T. Godoy, U. Kaufmann, \textit{On strictly positive solutions
for some semilinear elliptic problems}, NoDEA Nonlinear Differ. Equ. Appl.
\textbf{20} (2013), 779-795.

\bibitem {ans}T. Godoy, U. Kaufmann, \textit{Existence of strictly positive
solutions for sublinear elliptic problems in bounded domains}, Adv. Nonlinear
Stud. \textbf{14} (2014), 353-359.

\bibitem {jesusultimo}J. Hern\'{a}ndez, F. Mancebo, J. Vega, \textit{On the
linearization of some singular, nonlinear elliptic problems and applications},
Ann. Inst. H. Poincar\'{e} Anal. Non Lin\'{e}aire \textbf{19} (2002), 777--813.

\bibitem {J}L. Jeanjean, \textit{Some continuation properties via minimax
arguments}, Electron. J. Differential Equations \textbf{2011} (2011), Paper
No. 48, 10 pp.

\bibitem {K}R. Kajikiya, \textit{Positive solutions of semilinear elliptic
equations with small perturbations}, Proc. Amer. Math. Soc. \textbf{141}
(2013), 1335-1342.

\bibitem {LGMMT13}J. L\'{o}pez-G\'{o}mez, M. Molina-Meyer, A. Tellini,
\textit{The uniqueness of the linearly stable positive solution for a class of
superlinear indefinite problems with nonhomogeneous boundary conditions}, J.
Differential Equations \textbf{255} (2013), 503--523.

\bibitem {PT}M. A. Pozio and A. Tesei, \textit{Support properties of solution
for a class of degenerate parabolic problems}, Comm. Partial Differential
Equations \textbf{12} (1987), 47-75.

\bibitem {RQU6}H. Ramos Quoirin, K. Umezu, \textit{An indefinite
concave-convex equation under a Neumann boundary condition I},
arXiv:1603.04940 (to appear in Israel J.\ Math.).

\bibitem {tru}N. Trudinger, \textit{Linear elliptic operators with measurable
coefficients}, Ann. Scuola Norm. Sup. Pisa \textbf{27} (1973), 265--308.
\end{thebibliography}
\end{document}